\input ./style/arxiv-general.cfg
\documentclass[aos,MSNbibl,dvips]{arximspdf}
\makeatletter
   \@ifpackageloaded{graphicx}{}{\usepackage{graphicx}}
\makeatother
\usepackage{mathbh}

\doi{10.1214/15-AOS1270B}% Updated by VTEXPTS2LaTeX.exe, 26.02.2015
\volume{43}
\issue{4}
\pubyear{2015}
\firstpage{1437}
\lastpage{1443}
\docsubty{FLA}
\referstodoi{10.1214/14-AOS1270}

\makeatletter
\newcommand{\eqref}[1]{(\ref{#1})}
\makeatother

\begin{document}
\begin{frontmatter}
\vspace*{12pt}
\title{Discussion of ``Frequentist coverage of adaptive nonparametric
Bayesian credible sets''}
\runtitle{Discussion}

\begin{aug}
% Corresponding author: Ismael Castillo - ismael.castillo@upmc.fr% Updated by VTEXPTS2LaTeX.exe, 27.02.2015 08:32
%Updated by VTEXPTS2LaTeX.exe, 26.02.2015 12:55
\author{\fnms{Isma\"{e}l}~\snm{Castillo}\corref{}\ead[label=e1]{ismael.castillo@math.cnrs.fr}}%,
%\author[]{\fnms{}~\snm{}\ead[label=]{}}
% \and
%\author[]{\fnms{}~\snm{}\ead[label=]{}}
\runauthor{I. Castillo}
\affiliation{Universit\'es Paris VI \& VII}
%\dedicated{}
\address{CNRS---Laboratoire Probabilit\'es\\
\quad et Mod\`eles Al\'eatoires\\
Universit\'es Paris VI \& VII \\
B\^{a}timent Sophie Germain \\
75205 Paris Cedex 13\\
France\\
\printead{e1}}
%\address[]{\\\printead{}}
\end{aug}

% HISTORY:
%
\received{\smonth{1} \syear{2015}}% Updated by VTEXPTS2LaTeX.exe,
%26.02.2015 12:55
%\revised{\smonth{} \syear{}}

% ABSTRACT
%\begin{abstract}
%\end{abstract}

% KEYWORDS
%\begin{keyword}[class=AMS]
%\kwd[Primary ]{}
%\kwd{}
%\kwd[; secondary ]{}
%\end{keyword}
%\begin{keyword}
%\kwd{}
%\end{keyword}
\end{frontmatter}

First I would like to congratulate the three authors for a very nice
paper. During a visit to Eindhoven in 2010, Botond Szab\'o and Harry
van Zanten mentioned the first steps of this work to me, which
concerned the understanding of certain empirical Bayes procedures in
the white noise $L^2$-setting. Since then, together with Aad van der
Vaart, they have broadened their original goals and have produced an
impressive and very interesting series of papers on the subject. The
present paper is indeed one aspect of a larger body of work, and we
will mention a few connections with these related papers below.

The authors start from the signal in white noise model, that after
projection in $L^2$ onto an appropriate basis, typically related to the
SVD of the operator $K$ of the inverse problem, is translated into a
sequence formulation. They choose a prior distribution that makes
coordinates independent:
\begin{equation}
\label{def-prior} \Pi_\alpha\sim\bigotimes_{i\ge1}
N\bigl(0,i^{-1-2\alpha}\bigr).
\end{equation}
If the true parameter belongs to a regularity space defined from a
decay of coefficients in the previous basis, the authors prove that
certain credible sets constructed from the posterior distribution
coupled with a (marginal-likelihood) empirical Bayes (EB) procedure for
$\alpha$ achieve excellent performance: they are honest confidence sets
with adaptive, optimal asymptotic diameter if one restricts to certain
classes of ``self-similar''-type true parameters. These are the first
results of this type in Bayesian nonparametrics.

%Clearly, the paper opens the door to many interesting questions.
We organize this discussion around two main themes:
\begin{longlist}[2.]
\item[1.] Priors for Bayesian credible sets.
\item[2.] Bayesian credible regions and simulations.
\end{longlist}

\section{Priors for Bayesian credible sets}
Several aspects of the prior scheme \eqref{def-prior} are investigated
by the authors in \cite{kvv11}, \cite{ksvv12} together with Bartek
Knapik and in \cite{svv14}. In \cite{kvv11}, a fixed regularity
parameter $\alpha$ is considered; in \cite{ksvv12}, adaptative
contraction rates are derived. In \cite{svv14}, the prior \eqref
{def-prior} is used for fixed $\alpha$ and the use of a different
empirical Bayes scheme is advocated.

\begin{longlist}[{}]
\item[\textsc{Related priors.}] Staying with priors defined on the SVD of $K$,
some other adaptation schemes have been considered recently.
%Let us mention two of them.
One is (see \cite{svv13})
\begin{equation}
\label{priorsca} \Pi_\tau\sim\bigotimes_{i\ge1}
N\bigl(0,\tau^2i^{-1-2\alpha}\bigr),\qquad \tau>0,
\end{equation}
and adaptation is made by empirical Bayes or full Bayes on $\tau$.
%%This
%has the inconvenient of slightly restricting the range of adaptation
%depending on $\alpha$, namely adaptive rates can be derived from
%\eqref{priorsca} in the range $(0,\alpha+1/2)$ and not beyond.

Another prior is obtained by setting, for a sequence $\{\lambda_i\}
_{i\ge
1}$ of positive nondecreasing real numbers,
\begin{equation}
\label{priorheat} \Pi_t \sim\bigotimes_{i\ge1}
N\bigl(0,e^{- \lambda_i t}\bigr),\qquad t>0.
\end{equation}
In the case where $K$ is the identity and, for example, $\lambda_i=i^2$,
this falls into the framework considered in \cite{gid14}, where a full
Bayes method is considered by putting a well-chosen hyperprior on $t$.

A natural question is whether the same construction as in the paper
with a slightly blown up $L^2$-ball and regularity estimated by
empirical or full Bayes would work the same for the priors \eqref
{priorsca} or \eqref{priorheat}, with self-similarity constraints
expressed in a similar way. One can conjecture that the answer is yes
and that one may study the empirical Bayes procedure from the explicit
form of the marginal likelihood.

\item[\textsc{Related priors and sharp rates.}] Rates of convergence for Bayes
procedures are sometimes shown to be optimal up to a slowly varying
factor in $n$, for instance, logarithmic. In some cases it is not so
clear whether such a logarithmic term should be present in the rate or
not. The present work points to interesting questions with this
respect, with connections to the related\vspace*{1pt} prior schemes \eqref
{priorsca}--\eqref{priorheat}.

For prior \eqref{priorsca}, it is shown in \cite{svv13} that the
minimax rate $n^{-\beta/(1+2\beta)}$ in $L^2$ over hyperrectangles is
achieved by the marginal-likelihood-empirical Bayes procedure. This
comes, however, to a cost: one should assume that the true regularity
$\beta$ of the signal satisfies $\beta<1/2+\alpha$, for $\alpha$ the
regularity parameter in \eqref{priorsca}, otherwise the (uniform) rate
can be shown to be suboptimal.

For prior \eqref{priorheat}, we obtained in \cite{gid14} the rate
$(\log n/n)^{\beta/(1+2\beta)}$ in $L^2$ over a class containing
hyperrectangles and for which the minimax rate is $n^{-\beta/(1+2\beta)}$,
so \textit{without} the log-term, thus showing the unavoidable loss of a
logarithmic factor when using prior \eqref{priorheat}.

In \cite{ksvv12}, the authors obtain an upper-bound rate for prior
\eqref{def-prior} in $L^2$ that contains a logarithmic factor. However,
Proposition 3.8 of the present paper shows that the radius of the
credible set is proportional to $n^{-\beta/(1+2\beta)}$, while
Theorem~3.6
implies coverage of the credible set for polished tail parameters.
Combining these results, one deduces that the posterior mean $\hat
\theta
_{n,\hat\alpha_n}$ verifies $\| \hat\theta_{n,\hat\alpha_n} -
\theta_0\|_2=O_P(
n^{-\beta/(1+2\beta)})$. This presumably implies that the posterior itself
converges at the minimax rate, without extra log-terms, if the true
$\theta_0$ has polished tails. One may conjecture that this is also
true without the polished tail assumption. If so, it would be
interesting to better understand what makes that priors \eqref
{def-prior}--\eqref{priorheat} behave differently.

%There is something I do not yet understand in the comments end page
%15. There seems to be a logarithmic discrepancy between the
%convergence rate of the EB method as in \cite{ksvv12} and the radius
%of the credible set. How can this happen without harming coverage ? Is
%it due to the blowing up ? Probably rather to the polished-tails
%assumption of Theorem 3.6 (?) But I would have thought from our paper
%with do-ge that
%the lower bound with extra log factor in the rate was attained for a
%'self-similar' function ?? Or maybe for this prior there is no loss of
%a log ?? \\
%

\item[\textsc{Different priors and conditions}.] The prior scheme \eqref
{def-prior} is, by definition, somewhat tied to the SVD of $K$. As this
type of basis may not be well-localized, this may cause some
difficulties if the goal is a result in terms of a different loss
function than $L^2$.

Also, smoothness classes for $f_0$ are defined in terms of this basis
and thus connected to $K$. This may not always correspond to natural
assumptions of the practical problem at hand; see, for instance, \cite
{donoho95}.
The same can probably be said about the polished tail or
self-similarity conditions. As they stand, they refer to coefficients
in the basis associated to $K$, which may not always be canonical.

For these reasons, it would probably be interesting for future works to
consider different types of priors. It is unclear whether in general a
direct analysis of the explicit expression of the likelihood (and
marginal likelihood for the EB approach) will be possible. It would
certainly be desirable, if possible, to develop some general
understanding of empirical Bayes methods.
On the other hand, it would also be interesting to develop indirect (or
qualitative) techniques, similar to those of the meta-theorem of \cite
{ggv00} for these problems. Although this may not be easy for inverse
problems, some recent work for these include \cite{kr13} and \cite
{knasal14}. Other recent results on functionals using arguments
allowing implicit expressions can be found in \cite{cr14} and~\cite{ic14}.

%%%%%%%%%%% EXTRA %%%%%%%%%%%%%%%%%
% The comment that `It has already been noted in the literature that
%rates of contraction of functionals, such as a function at a point,
%are suboptimal unless the prior is made dependent on the functional'

\item[\textsc{Different approaches to nonparametric credible sets.}] As the
authors mention at the end of their introduction, for parametric models
the Bernstein--von Mises (BvM) theorem is a canonical tool to justify
that Bayesian credible sets are frequentist confidence sets. In \cite
{cn13} and \cite{cn14}, R. Nickl and myself proposed a possible
approach for the nonparametric BvM and showed that it could be applied
to the construction of fixed-regularity nonparametric confidence sets.
I am not sure I understand the authors' sentence
``no method that avoids dealing with the bias--variance trade-off will
properly quantify the uncertainty$\ldots$ current practice.'' In \cite
{cn13} and \cite{cn14}, no adaptation claims were made, and the
confidence sets there are for fixed regularity, although the proposed
methodology to build such sets does not  per se exclude adaptive priors.
Recently, a first application of this programme with adaptive priors in
white noise was carried out in \cite{rayprep}, leading to $L^2$ and
$L^\infty$ adaptive confidence sets computable in practice, under
appropriate self-similarity conditions. The ``bias--variance'' trade-off
mentioned by the authors I guess typically appears
when estimating the ``regularity'' of the signal, for instance, by an
empirical Bayes technique.

\item[\textsc{Bias--variance trade-off and choice of the prior}.] There are
several interesting questions mentioned by the authors beyond the
$L^2$-results of the paper. One is obtaining Bayesian confidence sets
for other norms, related to the problem of estimating certain
functionals, such as the value of the function at a point; see the
discussion on these in \cite{ksvv12} for the prior \eqref{def-prior}.
Another question is building different types of adaptive
$L^2$-confidence sets, where the regularity is assumed to belong to an
interval $[\alpha,2\alpha]$, as considered in \cite{svv14}, again
with the
scheme \eqref{def-prior}.

In both cases the
authors seem to conclude that marginal likelihood empirical Bayes or
full Bayes methods have some trouble, related to the choice of the
regularity parameter: for instance, the marginal-likelihood EB method
does not seem to perform the correct bias--variance trade-off in the two
problems. The proposed solution is then to choose the tuning parameter
$\hat\alpha_n$ independently, by a possibly non-Bayes method.
We agree, but one may note that all these results are for the given
prior scheme \eqref{def-prior}. Is it not conceivable that, for a given
problem (e.g., adaptive estimation of a functional), there exists a
prior for which the two steps are performed optimally? Perhaps this is
too much to ask in general, but, after all, this is the remarkable
result that the authors show in the present paper: at least for the
present problem, the Bayes method performs well in (1) rate-adaptation
and (2)~providing an (EB-)estimate $\hat\alpha_n$ so that the confidence
set has the desired coverage.
\end{longlist}

\section{Bayesian credible sets and simulations}

%{\sc Simulations.}
The authors present interesting simulations and a representation of the
credible sets in the case of the Volterra operator.

\begin{longlist}[{}]
\item[\textsc{What is exactly a plot of a credible set?}]
The credible ball considered in the paper is, with $L=1$,
\begin{equation}
\label{cs} \hat C_n = \bigl\{\theta\in\ell^2, \|
\theta-\hat\theta_{n,\hat\alpha
_n}\|_2\le r_{n,\gamma
}(\hat
\alpha_n) \bigr\}.
\end{equation}
In their Figure~1, the authors plot random draws from the posterior
distribution. The idea is that all (but possibly a few) of these draws
belong to the credible ball.
From this definition, we can make two comments:
\begin{longlist}[2.]
\item[1.] Curves that are not typical posterior draws belong to $\hat
C_n$.
\item[2.] There is typically much more ``information'' in the posterior
(coming from the prior) than the fact of belonging to such an $\ell^2$-ball.
\end{longlist}
To illustrate the fact that $\hat C_n$ is in some sense larger than the
``support'' of the posterior distribution, we have generated random
draws within $\hat C_n$ using a distribution different from the
posterior. First, consider the sequence, given the data,
\[
\mu=(\mu_k)_{k\ge1} \sim \biggl(\hat\theta_{n,\hat\alpha_n}+a
\frac{\xi_k}{(k\log^2
k)^{1/2}} \biggr)_{k\ge1},
\]
for $a>0$ some small constant and $\xi_k$ i.i.d. $N(0,1)$ variables.
Consider the law
\begin{equation}
\label{lawmu} \mathcal{L}(\mu|\mu\in\hat C_n),
\end{equation}
the distribution of $\mu$ conditioned to belong to the set $\hat C_n$.
Curves whose coefficients are sampled from this law are represented in
the left column of Figure~\ref{Figure1}, where we took $a=r_{n,\gamma
}(\hat\alpha_n)$, while the right column corresponds to posterior draws.
One notices that the typical curves on the left are more ``wiggly''
than those from the posterior distribution and also tend to spread
more, depending on how much curves $N$ are simulated, here $N=50$.

\begin{figure}

\includegraphics{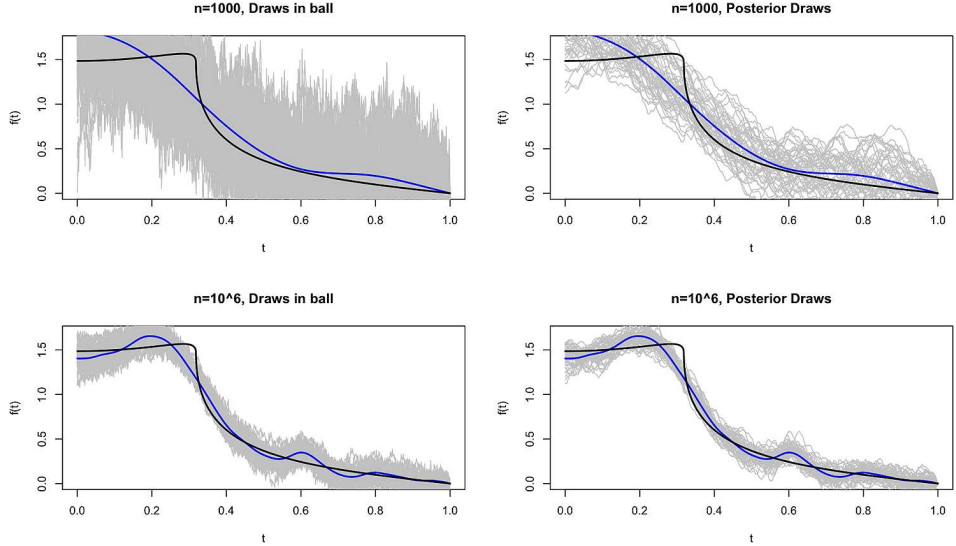}

\caption{In gray, on each plot, $N=50$ sampled curves from the
posterior distribution (right column) and the law induced by (\protect\ref{lawmu}) (left
column), for $n=10^3$ (top) and $n=10^6$ (bottom). Posterior mean and
true function are in blue and black,\vspace*{-3pt} respectively.}
\label{Figure1}
\end{figure}

On the other hand, the posterior distribution itself admits a series of
features that are not necessarily present
in a typical element of the $L^2$-ball. For instance, if $f$ is a draw
from the posterior on the signal function, and is $\hat\alpha_n$
concentrates, which is the case for self-similar-type truths, the
supremum norm $\|f-\hat f_{n,\hat\alpha_n}\|_\infty$ is a stochastically
bounded quantity that only depends on the data via $\hat\alpha_n$, as can
be seen from equation \eqref{law-rec} below. So with high probability
the posterior draws stay within a tube centered at the posterior mean.
If $\alpha>1/2$, one could presumably also prove at least some
supremum-norm consistency of the posterior around $f_0$, following, for
example, \cite{cr14}.

Given that the mathematical definition of the credible set is \eqref
{cs}, it seems natural to ask whether one should report draws from
posterior or from \eqref{lawmu}. Or rather, would it be possible to
define a credible set directly from the posterior draws themselves,
instead of reporting a full $L^2$-ball, while still retaining the
desired coverage properties?

%In theory one should report the L2 ball, in practice easy to report
%draws Can we do something= a draw of a plot in between that will still
%be fully justified theoretically (otherwise bunch of curves ...)

%Maybe one could 'bootstrap' some curves within the real 'credible' set
% based on $L^2$-constraint that the distance to the posterior mean
%should not be too large.
%
%The picture of random draws from the posterior actually produces some
%constraints which
%create a strict subset of the $L^2$ ball: some $L^\infty$ constraints
%in particular.
%{\em It seems that even the $L^2$ ball could have some $L^\infty$
%constraint in it depending on the basis (?) }

\item[\textsc{Improving on the estimate of the radius}.] The authors simulate\break 
$N=2000$ draws from the empirical Bayes posterior and retain the
$1-\gamma=95\%$ closest to the posterior mean. This means that an
implicit ``built-in'' estimator of the radius of the credible set is
used. More precisely, if $R_1,\ldots,R_N$ denote the observed
$L^2$-radii of $N$ draws under the posterior $\Pi_{\hat\alpha
_n}[\cdot
|X]$, only the curves with radius, respectively, $R_{(1)}\le\cdots
\le R_{(\lfloor 0.95\cdot N\rfloor)}$ are retained. In other words,
$R_{(\lfloor 0.95\cdot N\rfloor)}$ is used as an estimator of
$r_{n,\gamma
}(\hat\alpha_n)$.

This methodology\vspace*{1pt} is simple and certainly reasonable for relatively
large $N$, the precision of the ``built-in'' quantile estimator being
of order $N^{-1/2}$. In case one likes to be precise about the
$(1-\gamma
)$-coverage or, in cases where the posterior can only be approximated,
if one wants to detect possible outliers, one may suggest an
improvement based on a separate estimation of $r_{n,\gamma}(\hat
\alpha_n)$.
First, one may note that, in general, the posterior distribution of the
radius could be more easily accessible (or sampling from it could
require less computing time) than the full posterior. In the considered
white noise model example, computing a precise approximation of
$r_{n,\gamma}(\hat\alpha_n)$ is simple, as the posterior distribution
re-centered at the posterior mean has distribution, if $\tau_n$ is the map
$\theta\to\theta-\hat\theta_{n,\hat\alpha_n}$,
\begin{equation}
\label{law-rec} \Pi_{\hat\alpha_n}[\cdot|X]\circ\tau_n^{-1}
\,\mathop{=}^{\mathcal{L}}\, \bigotimes_{i\ge1} N
\biggl(0,\frac{1}{i^{1+2\hat\alpha
_n}+n\kappa
_i^2} \biggr).
\end{equation}
It is then straightforward to simulate the random variable $\|\zeta\|
_2$, where $\zeta$ is a draw from the distribution in the last display,
and then estimate $r_{n,\gamma}(\hat\alpha_n)$ based, for example,
on a
quantile as before, but this time using a much larger sample size (not
necessarily $N=2000$ as before). This can be made before running the
program simulating the posterior draws of the function $f$. For
instance, in the Volterra example with $n=1000$, one obtains the
estimate $\bar{r}_{n,\gamma}(1):=0.42\approx r_{n,\gamma}(1)$ using
a sample
of size $10^5$ [we set $\alpha=1$ for simplicity, but an approximation
of $r_{n,\gamma}(\hat\alpha_n)$ is obtained similarly, as soon as
$\hat\alpha_n$
has been computed].

We have run a few iterations of the algorithm proposed by the authors,
with the previous slight modification and setting $\alpha=1$ for
simplicity. As the estimate of the radius is improved, the rule for
discarding draws is more precise. For the results in Table~\ref{tata},
we have taken the precise estimate $\bar{r}_{n,\gamma}(1)$ as ``true.''

\begin{table}
\caption{Experiment using the original algorithm compared to a program
with separate precise estimation $\bar r_{n,\gamma}(1)$ (taken as
``true'') of $r_{n,\gamma}(1)$. After $10$ repetitions, $N_{fp}$ is the
mean number of ``incorrectly'' retained curves (false positive) by
original algorithm and $N_{fn}$ of ``incorrectly'' discarded curves
(false negative). In parenthesis percentage of occurrence}
\label{tata}
\begin{tabular*}{\tablewidth}{@{\extracolsep{\fill}}lcccccc@{}}
\hline
$\bolds{n}$ & \multicolumn{2}{c}{$\mathbf{1000}$} & \multicolumn{2}{c}{$\mathbf{10}^{\mathbf{6}}$} &
\multicolumn{2}{c@{}}{$\mathbf{10}^{\mathbf{8}}$} \\[-4pt]
& \multicolumn{2}{@{}l}{\hrulefill}  & \multicolumn{2}{l}{\hrulefill} & \multicolumn{2}{l@{}}{\hrulefill}\\
$\bolds{N}$ & \textbf{500} & \textbf{2000} & \textbf{500} & \textbf{2000} & \textbf{500} & \textbf{2000}\\
\hline
\\
$N_{fp}$ & $6$ (40\%) & \phantom{0}5 (70\%) & 4 $(50\%)$ & \phantom{0}6 $(40\%)$& 6 $(50\%)$
& 4 (20\%)\\
$N_{fn}$ & $3$ (50\%)  & 14 (20\%) &  6 (50\%)  & 12 (50\%)&  3 (50\%) & 8 (80\%)\\
\hline
\end{tabular*}
\end{table}
As shown in Table~\ref{tata}, a few curves per experiment typically
were either incorrectly included or excluded. Quantitatively, the
number of such curves is not very high,
but, on the other hand, these are the curves the farthest away from the
posterior mean, so visually this has (sometimes) some impact on the
pictures. This observation can be applied as well for pictures of
credible bands, as recently considered, for example, in \cite{rayprep}.

%The probability that there is at least a draw in this picture that
%either A) is actually not in the credible set
%as defined or that B) we have discarded a draw that actually belongs
%to it, is probably quite high. For A) it is possibly visually not too
%good as `bad' draws are probably the ones one sees most. B) could be a
%bit annoying too, as the plotted set could possibly be slightly
%overconfident by missing some `not-too-good' draws.
%(As $n$ grows, both tend to $1/2$, as $\mathbb{P}(
%\text{Binomial}(n,.05)=n*.05)$ goes to $0$ with $n$ and $\mathbb{P}(
%\text{Binomial}(n,.05)>n*.05)$ is close to $1/2$ -check median-)
%But here we actually {\em know} (at least possibly with extremely high
%precision) which of the draws {\em actually belong} to the credible
%set. Indeed, to do that it is enough to
%simulate from the (recentered) $L^2$ radius of the posterior. Here the
%law is the $L^2$ norm of a Gaussian measure with given parameters
%(depending on $\hat\alpha$ which, given the data, is fixed), so the
%law could be tabulated.
%-check whether there are some classical references on quantiles of
%$L^2$ norms of Gaussian vectors or a way to approximate them
%efficiently. The Monte Carlo error could be made quite low I guess. \\

Congratulations again to the authors for their inspiring series of
works. Developing tools to build Bayesian credible sets for other
models and priors is a very interesting topic, and we expect to see
more on the subject soon.
%Developing further insight on prior families appears to be important.
%
%From this very nice body of work one may ask a few questions for
%future work
%\begin{enumerate}
%\item Qualitative insight on different empirical Bayes procedures
%\item Provide tools for estimators for which direct computation is
%difficult
%\item Derive results in other models or for other norms (bands etc.)
%\end{enumerate}
\end{longlist}

% imsref loaded by daiva.urboniene, 2015-02-26 13:18:30

%\begin{appendix}
%\section{}
%\end{appendix}

%\section*{Acknowledgments}

%\begin{supplement}[id=suppA]
%\slink[doi]{10.1214/00-AOSXXXXSUPP} %[doi,text={...}] - jei reikia
%suskaldyti doi
%\sfilename{aosXXXX\_supp.pdf}
%\end{supplement}

%\begin{thebibliography}{99}
%\bibitem[\protect\citeauthoryear{}{}]{r1}
%\bibitem{r1}
%\end{thebibliography}

\printaddresses
\end{document}